\documentclass[12pt]{amsart}

\usepackage{enumitem}
\usepackage{psfrag}
\usepackage{graphicx}
\usepackage{pinlabel}
\usepackage{graphicx}
\usepackage{amsmath}
\usepackage{amssymb}
\usepackage{amscd}
\usepackage{autobreak}
\usepackage{picinpar}
\usepackage{times}
\usepackage{pb-diagram}
\usepackage{graphicx}
\usepackage{wrapfig}
\usepackage{xspace}
\usepackage{hyperref}
\usepackage{url}
\usepackage{caption}
\usepackage{subcaption}
\usepackage{fancyhdr}
\usepackage{overpic}
\usepackage{color}
\usepackage[square,comma,sort&compress,numbers]{natbib}
\usepackage{latexsym}
\usepackage{mathrsfs}
\usepackage{amsfonts}
\usepackage{hyperref}
\usepackage{amsthm}
\usepackage{appendix}
\usepackage{pdfsync}

\theoremstyle{definition}
\newtheorem{theorem}{Theorem}[section]

\newtheorem{definition}[theorem]{Definition}

\newtheorem{remark}[theorem]{Remark}
\newtheorem*{theorem*}{Theorem}
\def\qed{\hfill{Q.E.D.}\smallskip}

\begin{document}

\title{\bf On the classification of discrete conformal structures on surfaces}
\author{Xu Xu, Chao Zheng}

\date{\today}

\address{School of Mathematics and Statistics, Wuhan University, Wuhan, 430072, P.R.China}
 \email{xuxu2@whu.edu.cn}

\address{School of Mathematics and Statistics, Wuhan University, Wuhan 430072, P.R. China}
\email{czheng@whu.edu.cn}


\begin{abstract}
Glickenstein \cite{Glickenstein} and Glickenstein-Thomas \cite{GT} introduced the discrete conformal structures on surfaces in an axiomatic approach
and studied its classification.
In this paper, we give a full classification of the discrete conformal structures on surfaces, which completes Glickenstein-Thomas' classification.
As a result, we find some new classes of discrete conformal structures on surfaces, including some of the generalized circle packing metrics introduced by Guo-Luo \cite{GL2}.
The relationships between the discrete conformal structures on surfaces and the 3-dimensional hyperbolic geometry are also discussed.
\end{abstract}

\maketitle

\section{Introduction}

Discrete conformal structure on polyhedral surfaces is a discrete analogue of the smooth conformal structure on surfaces that assigns discrete metrics by scalar functions defined on the vertices.
The generic notion of discrete conformal structure on polyhedral surfaces was introduced independently by Glickenstein \cite{Glickenstein} and Glickenstein-Thomas \cite{GT} in an axiomatic approach.
Glickenstein \cite{Glickenstein} first introduced
the notion of partial edge length for a closed triangulated surface in the Euclidean background geometry, which ensures the existence of some geometric structures on the Poincar\'{e} dual of the triangulated surface.
Following \cite{Glickenstein}, Glickenstein-Thomas \cite{GT} further introduced the partial edge length for a closed triangulated surface in the hyperbolic and spherical background geometry.
Based on the partial edge length,
Glickenstein \cite{Glickenstein} and Glickenstein-Thomas \cite{GT} further introduced the definition of discrete conformal structures in an axiomatic approach, in which the deformation of discrete conformal structure is required to depend in a reasonable form of the partial edge length.
In  \cite{GT}, Glickenstein-Thomas gave a classification of discrete conformal structures on surfaces.
However, Glickenstein-Thomas's classification of  the discrete conformal structures on surfaces in the hyperbolic and spherical background geometry
is not complete.
In this paper, we give a full classification of the discrete conformal structures on surfaces, which
completes Glickenstein-Thomas's classification.

Suppose $(M, \mathcal{T})$ is a triangulated connected closed surface with a triangulation $\mathcal{T}={(V,E,F)}$, where $V,E,F$ represent the sets of vertices, edges and faces respectively and $V$ is a finite subset of $M$.
Denote a vertex, an edge and a face in the triangulation $\mathcal{T}$ by $v_i, \{ij\}, \{ijk\}$ respectively, where ${i,j,k}$ are natural numbers.
We use $E_{+}$ to denote oriented edges and label them as ordered pairs $(i,j)\in E_{+}$.
The sets of real valued functions on $V, E$ and $E_{+}$ are denoted by $V^*$, $E^*$ and $E_+^*$ respectively.
If a map $l: E\rightarrow (0, +\infty)$ assigns a length to every edge in such a way that every triangle $\{ijk\}\in F$ with edge lengths $l_{ij}, l_{ik}, l_{jk}$ is embedded in $\mathbb{G}$ ($\mathbb{G}=\mathbb{E}^2, \mathbb{H}^2$ or $\mathbb{S}^2$),
then $(M, \mathcal{T}, l)$ is called as a piececwise $\mathbb{G}$ surface and $l: E\rightarrow (0, +\infty)$ is called a $\mathbb{G}$ polyhedral metric.
One can also take the triangulated surface $(M, \mathcal{T}, l)$ being obtained by gluing triangles in $\mathbb{E}^2$ ($\mathbb{H}^2$ or $\mathbb{S}^2$ respectively) isometrically along the edges in pair.

To ensure the existence of a dual geometric structure on the Poincar\'{e} dual of a triangulated surface $(M,\mathcal{T})$ with $\mathbb{G}$ background geometry,
Glickenstein \cite{Glickenstein} and Glickenstein-Thomas \cite{GT} introduced the following definition of partial edge length for a triangulated surface.

\begin{definition}[\cite{Glickenstein,GT}]\label{Def: partial edge length}
Suppose $(M,\mathcal{T})$ is a triangulated surface.
An element $d\in E_+^*$ is an assignment of partial edge lengths, if $l_{ij}=d_{ij}+d_{ji}$ makes $(M, \mathcal{T}, l)$ a piecewise $\mathbb{G}$ surface and if for every triangle $\{ijk\}\in F$,
\begin{equation*}
d^2_{ij}+d^2_{jk}+d^2_{ki}=d^2_{ji}+d^2_{kj}+d^2_{ik}   \quad  \text{if $\mathbb{G}=\mathbb{E}^2$};
\end{equation*}
\begin{equation}\label{Eq: hyperbolic compatible condition}
\cosh(d_{ij})\cosh(d_{jk})\cosh(d_{ki})
=\cosh(d_{ji})\cosh(d_{kj})\cosh(d_{ik})   \quad  \text{if $\mathbb{G}=\mathbb{H}^2$};
\end{equation}
\begin{equation}\label{Eq: spherical compatible condition}
\cos(d_{ij})\cos(d_{jk})\cos(d_{ki})
=\cos(d_{ji})\cos(d_{kj})\cos(d_{ik})   \quad  \text{if $\mathbb{G}=\mathbb{S}^2$}.
\end{equation}
%We say that $(M,\mathcal{T},d)$ is a $\mathbb{G}-$metrized surface.
\end{definition}

One can define the edge center $c_{ij}$ on the geodesic $E_{ij}$ extending each edge $\{ij\}\in E$ of a triangle $\{ijk\}\in F$ with edge lengths $l_{ij}, l_{ik}, l_{jk}$ in $\mathbb{G}$ background geometry, which is the unique point that is of signed distance $d_{ij}$ to $v_i\in V$ and $d_{ji}$ to $v_j\in V$.
Note that the signed distance $d_{ij}$ of $c_{ij}$ to  $v_i\in V$ is positive if $c_{ij}$ is on the same side of $v_i$ as $v_j$ along the geodesic $E_{ij}$, and is negative otherwise.

A discrete conformal structure on a triangulated surface $(M,\mathcal{T})$ assigns partial edge length $d$ defined on the oriented edges $E_+$ by scalar functions defined on the vertices $V$.
Glickenstein \cite{Glickenstein} and  Glickenstein-Thomas \cite{GT} introduced the following definition of discrete conformal structures on surfaces in an axiomatic approach.
\begin{definition}[\cite{Glickenstein,GT}]\label{Def: DCS}
Suppose $(M,\mathcal{T},l)$ is a piecewise $\mathbb{G}$ surface.
A discrete conformal structure $d=d(f)$ on $(M,\mathcal{T},l)$ is a smooth map from $U\subseteq V^*$ to $E_{+}^*$ that sends functions $f$ defined on the vertices to a set of partial edge lengths $\{d_{ij}\}$ such that for each $(i,j)\in E_{+}$ and $k\in V$,
\begin{equation*}%\label{Eq: Euclidean DCS}
\frac{\partial l_{ij}}{\partial f_i}=d_{ij} \quad  \text{if $\mathbb{G}=\mathbb{E}^2$},
\end{equation*}
\begin{equation}\label{Eq: hyperbolic DCS}
\frac{\partial l_{ij}}{\partial f_i}=\tanh d_{ij} \quad  \text{if $\mathbb{G}=\mathbb{H}^2$},
\end{equation}
\begin{equation}\label{Eq: spherical DCS}
\frac{\partial l_{ij}}{\partial f_i}=\tan d_{ij} \quad  \text{if $\mathbb{G}=\mathbb{S}^2$},
\end{equation}
and
\begin{equation}\label{Eq: variation}
\frac{\partial d_{ij}}{\partial f_k}=0
\end{equation}
if $k\neq i$ and $k\neq j$.
\end{definition}

The main result of this paper is the following classification theorem, which completes Glickenstein-Thomas' classification of the discrete conformal structures on surfaces \cite{GT}.

\begin{theorem}\label{Thm: classification of DCS}
Let $d=d(f)$ be a discrete conformal structure on a triangulated surface $(M,\mathcal{T})$ with $\mathbb{G}$ background geometry.
Then there exist $\alpha\in V^*$ and $\eta\in E^*$ with $\eta_{ij}=\eta_{ji}$ such that
\begin{description}
\item[(A)] if $\mathbb{G}=\mathbb{E}^2$, then the discrete conformal structure $d=d(f)$ can be written as
\begin{equation*}
d_{ij}=\frac{\alpha_ie^{2f_i}+\eta_{ij}e^{f_i+f_j}}{l_{ij}}
\end{equation*}
with
\begin{equation*}
l^2_{ij}=\alpha_ie^{2f_i}+\alpha_je^{2f_j}+2\eta_{ij}e^{f_i+f_j};
\end{equation*}

\item[(B)] if $\mathbb{G}=\mathbb{H}^2$, then the discrete conformal structure $d=d(f)$  has one of the following two forms,
\begin{description}
\item[(b1)]
\begin{equation*}
\tanh d_{ij}=\frac{\alpha_ie^{2f_i}}{\sinh l_{ij}}\sqrt{\frac{1+\alpha_je^{2f_j}}{1+\alpha_ie^{2f_i}}}
+\frac{\eta_{ij}e^{f_i+f_j}}{\sinh l_{ij}}
\end{equation*}
with
\begin{equation}\label{Eq: F22}
\cosh l_{ij}
=\sqrt{(1+\alpha_ie^{2f_i})(1+\alpha_je^{2f_j})}
+\eta_{ij}e^{f_i+f_j},\ \text{for}\ 1+\alpha_ie^{2f_i}>0, 1+\alpha_je^{2f_j}>0,
\end{equation}
\begin{equation}\label{Eq: F23}
\cosh l_{ij}
=-\sqrt{(1+\alpha_ie^{2f_i})(1+\alpha_je^{2f_j})}
+\eta_{ij}e^{f_i+f_j},\ \text{for}\ 1+\alpha_ie^{2f_i}<0, 1+\alpha_je^{2f_j}<0.
\end{equation}

\item[(b2)]
\begin{equation*}
\tanh d_{ij}=-\frac{\sinh(f_j-f_i-C_{ij})}{\sinh l_{ij}}+\frac{\eta_{ij}e^{f_i+f_j}}{\sinh l_{ij}}
\end{equation*}
with
\begin{equation}\label{Eq: F25}
\cosh l_{ij}
=\cosh(f_j-f_i-C_{ij})+\eta_{ij}e^{f_i+f_j},
\end{equation}
where $C_{ij}$ is a constant satisfying $C_{ij}+C_{jk}+C_{ki}=0$ and $C_{rs}+C_{sr}=0$ for $\{r,s\}\subseteq\{i,j,k\}$;
\end{description}

\item[(C)] if $\mathbb{G}=\mathbb{S}^2$, then
when $\frac{\cos d_{ij}}{\cos d_{ji}}>0$, the discrete conformal structure $d=d(f)$  has one of the following two forms,
\begin{description}
\item[(c1)]
\begin{equation*}
\tan d_{ij}=\frac{\alpha_ie^{2f_i}}{\sin l_{ij}}\sqrt{\frac{1-\alpha_je^{2f_j}}{1-\alpha_ie^{2f_i}}}
+\frac{\eta_{ij}e^{f_i+f_j}}{\sin l_{ij}}
\end{equation*}
with
\begin{equation}\label{Eq: F32}
\cos l_{ij}
=\sqrt{(1-\alpha_ie^{2f_i})(1-\alpha_je^{2f_j})}
-\eta_{ij}e^{f_i+f_j},\ \text{for}\ 1-\alpha_ie^{2f_i}>0, 1-\alpha_je^{2f_j}>0,
\end{equation}
\begin{equation}\label{Eq: F33}
\cos l_{ij}
=-\sqrt{(1-\alpha_ie^{2f_i})(1-\alpha_je^{2f_j})}
-\eta_{ij}e^{f_i+f_j},\ \text{for}\ 1-\alpha_ie^{2f_i}<0, 1-\alpha_je^{2f_j}<0.
\end{equation}

\item[(c2)]
\begin{equation*}
\tan d_{ij}=\frac{\sinh(f_j-f_i-C_{ij})}{\sin l_{ij}}+\frac{\eta_{ij}e^{f_i+f_j}}{\sin l_{ij}}
\end{equation*}
with
\begin{equation*}
\cos l_{ij}
=\cosh(f_j-f_i-C_{ij})-\eta_{ij}e^{f_i+f_j},
\end{equation*}
where $C_{ij}$ is a constant satisfying $C_{ij}+C_{jk}+C_{ki}=0$ and $C_{rs}+C_{sr}=0$ for $\{r,s\}\subseteq\{i,j,k\}$;
\end{description}
when $\frac{\cos d_{ij}}{\cos d_{ji}}<0$, the discrete conformal structure $d=d(f)$  has one of the following two forms,
\begin{description}
\item[(c3)]
\begin{equation*}
\tan d_{ij}=-\frac{\alpha_ie^{2f_i}}{\sin l_{ij}}\sqrt{\frac{1-\alpha_je^{2f_j}}{1-\alpha_ie^{2f_i}}}
-\frac{\eta_{ij}e^{f_i+f_j}}{\sin l_{ij}}
\end{equation*}
with
\begin{equation}\label{Eq: F34}
\cos l_{ij}
=-\sqrt{(1-\alpha_ie^{2f_i})(1-\alpha_je^{2f_j})}
+\eta_{ij}e^{f_i+f_j},\ \text{for}\ 1-\alpha_ie^{2f_i}>0, 1-\alpha_je^{2f_j}>0,
\end{equation}
\begin{equation}\label{Eq: F35}
\cos l_{ij}
=\sqrt{(1-\alpha_ie^{2f_i})(1-\alpha_je^{2f_j})}
+\eta_{ij}e^{f_i+f_j},\ \text{for}\ 1-\alpha_ie^{2f_i}<0, 1-\alpha_je^{2f_j}<0.
\end{equation}

\item[(c4)]
\begin{equation*}
\tan d_{ij}=-\frac{\sinh(f_j-f_i-C_{ij})}{\sin l_{ij}}-\frac{\eta_{ij}e^{f_i+f_j}}{\sin l_{ij}}
\end{equation*}
with
\begin{equation*}
\cos l_{ij}
=-\cosh(f_j-f_i-C_{ij})+\eta_{ij}e^{f_i+f_j},
\end{equation*}
where $C_{ij}$ is a constant satisfying $C_{ij}+C_{jk}+C_{ki}=0$ and $C_{rs}+C_{sr}=0$ for $\{r,s\}\subseteq\{i,j,k\}$.
\end{description}
\end{description}
\end{theorem}

\begin{remark}
Glickenstein-Thomas' classification theorem in \cite{GT} contains only the discrete conformal structures $\mathbf{(A)}$, (\ref{Eq: F22}) and (\ref{Eq: F32}) in Theorem \ref{Thm: classification of DCS}.
We add the discrete conformal structures (\ref{Eq: F23}), $\mathbf{(b2)}$, (\ref{Eq: F33}), $\mathbf{(c2)}$, $\mathbf{(c3)}$ and $\mathbf{(c4)}$ in Theorem \ref{Thm: classification of DCS} to complete the classification of discrete conformal structure on a triangulated surface with $\mathbb{G}$ background geometry.
As the proof for the discrete conformal structures $\mathbf{(A)}$ has been given in details in \cite{GT, Thomas}, so we only give the proof for the discrete conformal structures $\mathbf{(B)}$ and $\mathbf{(C)}$ in this paper.
\end{remark}

\begin{remark}\label{remark: 2}
As pointed out by Thomas (\cite{Thomas}, page 53),
one can reparameterize discrete conformal structures $d=d(f)$ in $\mathbf{(A)}$, (\ref{Eq: F22}) and (\ref{Eq: F32}) of Theorem \ref{Thm: classification of DCS} so that $\alpha_i\in \{-1,0,1\}$ while keeping the induced $\mathbb{G}$ polyhedral metrics invariant.
Similar parameterization also applies to the discrete conformal structure (\ref{Eq: F34}).
The conditions that $1+\alpha_ie^{2f_i}<0, 1+\alpha_je^{2f_j}<0$ in the discrete conformal structures (\ref{Eq: F23}) imply $\alpha_i<0, \alpha_j<0$.
Thus we always assume $\alpha\equiv -1$ in the discrete conformal structures (\ref{Eq: F23}) in the following.
Similarly, we always assume $\alpha\equiv +1$ in the discrete conformal structures (\ref{Eq: F33}) and (\ref{Eq: F35}) in the following.

Moreover, if a triangulated surface $(M,\mathcal{T})$ with $\mathbb{G}\ (=\mathbb{H}^2\ \text{or}\ \mathbb{S}^2)$ background geometry is simply connected,
one can also reparameterize discrete conformal structures $d=d(f)$ in $\mathbf{(b2)}$, $\mathbf{(c2)}$ and $\mathbf{(c4)}$ so that $C_{ij}=0$ while keeping the induced $\mathbb{G}$ polyhedral metrics invariant respectively.
Indeed, for these discrete conformal structures in  $\mathbf{(b2)}$, $\mathbf{(c2)}$ and $\mathbf{(c4)}$, we have
\begin{eqnarray}\label{Eq: F21}
\begin{cases}
C_{ij}+C_{jk}+C_{ki}=0,\\
C_{ij}+C_{ji}=0,\\
C_{jk}+C_{kj}=0,\\
C_{ki}+C_{ik}=0
\end{cases}
\end{eqnarray}
for any triangle $\{ijk\}\in F$.
Fix a vertex $v_0$.
By the connectedness of $(M, \mathcal{T})$, for any vertex $v\in V$, there exists a path $\gamma_1: v_0\sim v_1\sim...\sim v_{n-1}\sim v_n=v$ connecting $v_0$ and $v$.
We define a function   $g: V\rightarrow \mathbb{R}$ by setting   $g(v)=g_n=-\sum_{i=1}^nC_{i-1,i}$ and $g_0=0$.
Then $g$ is well-defined.
In fact, if there exists another path $\gamma_2: v_0\sim v^\prime_1\sim...\sim v^\prime_{m-1}\sim v^\prime_m=v$ connecting $v_0$ and $v$,
then $\widetilde{g}_m=-\sum_{k=1}^mC^\prime_{k-1,k}$ by definition.
By the simply connectedness of $(M,\mathcal{T})$ and the condition (\ref{Eq: F21}), we have  $g_n=\widetilde{g}_m$. Therefore, $C_{ij}=g_{i}-g_j$ for $i$ adjacent to $j$.
We claim $g$ is unique up to a constant.
Suppose there exists another function $\overline{g}: V\rightarrow \mathbb{R}$ such that
$C_{ij}=\overline{g}_{i}-\overline{g}_j$ for $i$ adjacent to $j$.
Then $g_{i}-\overline{g}_{i}=g_j-\overline{g}_j$  for $i$ adjacent to $j$.
By the connectedness of $(M,\mathcal{T})$, we have $g-\overline{g}=c(1,...,1)$, where $c$ is a constant.
Therefore, (\ref{Eq: F25}) is equivalent to
\begin{equation*}
\cosh l_{ij}
=\cosh(f_j-f_i-g_i+g_j)+\eta_{ij}e^{f_i+f_j}
=\cosh(h_j-h_i)+\widetilde{\eta}_{ij}e^{h_i+h_j},
\end{equation*}
where $h_i=f_i+g_i$, $\widetilde{\eta}_{ij}=e^{-g_i-g_j}\eta_{ij}$.
Similar arguments apply to the discrete conformal structures $\mathbf{(c2)}$ and $\mathbf{(c4)}$.
\end{remark}

\begin{remark}\label{remark 1}
It should be mentioned that the discrete conformal structures $\mathbf{(b1)}$ and $\mathbf{(b2)}$ can not exist simultaneously on the same triangulated surface $(M,\mathcal{T})$ by the proof of Theorem \ref{Thm: classification of DCS} $\mathbf{(B)}$ in Section \ref{section 2}.
And the discrete conformal structures $\mathbf{(c1)}$ and $\mathbf{(c2)}$ can not exist simultaneously on the same triangulated surface $(M,\mathcal{T})$ by the proof of Theorem \ref{Thm: classification of DCS} $\mathbf{(C)}$ in Section \ref{section 3}.
It is possible that the discrete conformal structures $\mathbf{(c1)}$ and $\mathbf{(c2)}$ exist alone on a triangulated surface $(M,\mathcal{T})$ with spherical background geometry.
According to Theorem \ref{Thm: classification of DCS} again, the discrete conformal structures $\mathbf{(c3)}$ and $\mathbf{(c4)}$ exist in $(M,\mathcal{T})$ only if $\frac{\cos d_{ij}}{\cos d_{ji}}<0$ for some spherical triangles $\{ijk\}\in F$.
This implies that the discrete conformal structures $\mathbf{(c3)}$ and $\mathbf{(c4)}$ can not exist alone on a triangulated surface $(M,\mathcal{T})$ with spherical background geometry.
Otherwise, we have $\frac{\cos d_{ij}}{\cos d_{ji}}<0$,\ $\frac{\cos d_{ik}}{\cos d_{ki}}<0$ and $\frac{\cos d_{jk}}{\cos d_{kj}}<0$, which contradicts the condition (\ref{Eq: spherical compatible condition}).
Therefore, we have to consider the mixed type of discrete conformal structures, i.e., the discrete hyperbolic metrics $l_{ij}, l_{jk}, l_{ki}$ of a spherical triangle $\{ijk\}\in F$ are induced by $\mathbf{(c1)}$ and $\mathbf{(c3)}$ or $\mathbf{(c2)}$ and $\mathbf{(c4)}$, such that the condition (\ref{Eq: spherical compatible condition}) holds.
\end{remark}

The paper is organized as follows.
In Section \ref{section 2}, we prove the hyperbolic case in Theorem \ref{Thm: classification of DCS}.
In Section \ref{section 3}, we prove the spherical case in Theorem \ref{Thm: classification of DCS}.
In Section \ref{section 4}, we discuss the relationships between discrete conformal structures in Theorem \ref{Thm: classification of DCS} and the 3-dimensional hyperbolic geometry.

\section{Proof of hyperbolic case in Theorem \ref{Thm: classification of DCS} $\mathbf{(B)}$}\label{section 2}

For the discrete conformal structure $\mathbf{(b1)}$, Glickenstein-Thomas \cite{GT} did not give a detailed proof, although the proof is similar.
For completeness, we provide the proof of the case $\mathbf{(B)}$ in details.
By (\ref{Eq: hyperbolic DCS}), Glickenstein-Thomas \cite{GT} obtained
\begin{equation}\label{Eq: F1}
\frac{\partial}{\partial f_i}\cosh l_{ij}
=\cosh l_{ij}-\frac{\cosh d_{ji}}{\cosh d_{ij}},
\end{equation}
\begin{equation}\label{Eq: F2}
\frac{\partial}{\partial f_j}\cosh l_{ij}
=\cosh l_{ij}-\frac{\cosh d_{ij}}{\cosh d_{ji}}.
\end{equation}
Furthermore, they gave
\begin{equation}\label{Eq: F3}
\frac{1}{\cosh^2 d_{ij}}\frac{\partial d_{ij}}{\partial f_j}=\frac{\partial}{\partial f_j}\tanh d_{ij}
=\frac{\partial^2 l_{ij}}{\partial f_i\partial f_j}
=\frac{\partial}{\partial f_i}\tanh d_{ji}=\frac{1}{\cosh^2 d_{ji}}\frac{\partial d_{ji}}{\partial f_i}.
\end{equation}
Set
\begin{equation}\label{defn of H}
H=\log\frac{\cosh^2 d_{ij}}{\cosh^2 d_{ji}}.
\end{equation}
Then the formulas (\ref{Eq: hyperbolic compatible condition}) and (\ref{Eq: variation}) give
\begin{equation}\label{Eq: F4}
\frac{\partial^2 H}{\partial f_i\partial f_j}=0.
\end{equation}
By direct but tedious calculations, Glickenstein-Thomas \cite{GT} obtained the following formula
\begin{equation}\label{Eq: F18}
\bigg(\frac{\cosh^2 d_{ij}}{\cosh^2 d_{ji}}\frac{\partial}{\partial f_i}+\frac{\partial}{\partial f_j}\bigg)H=2\bigg(\frac{\cosh^2 d_{ij}}{\cosh^2 d_{ji}}-1\bigg).
\end{equation}
One can refer to Thomas \cite{Thomas} for the details of the calculations.
The formula (\ref{Eq: F18}) is equivalent to
\begin{equation}\label{Eq: F5}
\bigg(e^H\frac{\partial}{\partial f_i}+\frac{\partial}{\partial f_j}\bigg)H=2(e^H-1),
\end{equation}
and
\begin{equation}\label{Eq: F6}
\bigg(\frac{\partial}{\partial f_i}+e^{-H}\frac{\partial}{\partial f_j}\bigg)H=2(1-e^{-H}).
\end{equation}
Furthermore, using (\ref{Eq: F4}) and differentiating (\ref{Eq: F5}) with respect to $f_i$ and (\ref{Eq: F6}) with respect to $f_j$, Glickenstein-Thomas \cite{GT} obtianed
\begin{equation}\label{Eq: partial H partial fi}
\frac{\partial^2 H}{\partial f^2_i}+\bigg(\frac{\partial H}{\partial f_i}\bigg)^2=2\frac{\partial H}{\partial f_i},
\end{equation}
and
\begin{equation}\label{Eq: partial H partial fj}
\frac{\partial^2 H}{\partial f^2_j}-\bigg(\frac{\partial H}{\partial f_j}\bigg)^2=2\frac{\partial H}{\partial f_j}.
\end{equation}
One can easily solve the ODE (\ref{Eq: partial H partial fi}) to get that
\begin{equation*}
\frac{\partial H}{\partial f_i}\equiv2\ \quad \text{or}\ \quad
\frac{\partial H}{\partial f_i}=2\frac{a_{ij}e^{2f_i}}{1+a_{ij}e^{2f_i}}
\end{equation*}
for some constants $a_{ij}$ by (\ref{Eq: F4}), and solve the ODE (\ref{Eq: partial H partial fj}) to get that
\begin{equation*}
\frac{\partial H}{\partial f_j}\equiv -2\ \quad \text{or}\ \quad
\frac{\partial H}{\partial f_j}=-2\frac{a_{ji}e^{2f_j}}{1+a_{ji}e^{2f_j}}
\end{equation*}
for some constants $a_{ji}$ by (\ref{Eq: F4}).
Note that the formula (\ref{Eq: F5}) can be rewritten as
\begin{equation}\label{Eq: F8}
(\frac{\partial H}{\partial f_i}-2)e^H
=-(\frac{\partial H}{\partial f_j}+2).
\end{equation}
\begin{description}
\item[(I)]
If $\frac{\partial H}{\partial f_i}\not\equiv 2$, then $\frac{\partial H}{\partial f_j}\not\equiv -2$ by (\ref{Eq: F8}).
This implies that $\frac{\partial H}{\partial f_i}=\frac{2a_{ij}e^{2f_i}}{1+a_{ij}e^{2f_i}}$ and $\frac{\partial H}{\partial f_j}=-\frac{2a_{ji}e^{2f_j}}{1+a_{ji}e^{2f_j}}$.
Then
\begin{equation*}
\frac{\cosh^2 d_{ij}}{\cosh^2 d_{ji}}
=e^H
=\frac{1}{\frac{\partial H}{\partial f_i}-2}(-\frac{\partial H}{\partial f_j}-2)
=\frac{1+a_{ij}e^{2f_i}}{1+a_{ji}e^{2f_j}},
\end{equation*}
Using (\ref{Eq: F1}) and (\ref{Eq: F2}), one can obtain the following system
\begin{equation*}
\begin{aligned}
&\frac{\partial}{\partial f_i}\cosh l_{ij}
=\cosh l_{ij}
-\sqrt{\frac{1+a_{ji}e^{2f_j}}{1+a_{ij}e^{2f_i}}},\\
&\frac{\partial}{\partial f_j}\cosh l_{ij}
=\cosh l_{ij}
-\sqrt{\frac{1+a_{ij}e^{2f_i}}{1+a_{ji}e^{2f_j}}}.
\end{aligned}
\end{equation*}
This implies
\begin{gather*}
\cosh l_{ij}
=\sqrt{(1+a_{ij}e^{2f_i})(1+a_{ji}e^{2f_j})}
+\eta_{ij}e^{f_i+f_j},\ \text{if}\ 1+a_{ij}e^{2f_i}>0,\ 1+a_{ji}e^{2f_j}>0, \\
\cosh l_{ij}
=-\sqrt{(1+a_{ij}e^{2f_i})(1+a_{ji}e^{2f_j})}
+\eta_{ij}e^{f_i+f_j},\ \text{if}\ 1+a_{ij}e^{2f_i}<0,\ 1+a_{ji}e^{2f_j}<0, \\
\end{gather*}
where $\eta_{ij}\in \mathbb{R}$ is a constant.
By (\ref{Eq: hyperbolic compatible condition}) and (\ref{Eq: F4}), we have
\begin{equation}\label{Eq: F31}
\log \frac{\cosh d_{ij}}{\cosh d_{ji}}
+\log \frac{\cosh d_{ki}}{\cosh d_{ik}}
=\log \frac{\cosh d_{kj}}{\cosh d_{jk}}.
\end{equation}
Note that $\log \frac{\cosh d_{kj}}{\cosh d_{jk}}$ is independent of $f_i$ by (\ref{Eq: variation}),
so differentiating (\ref{Eq: F31}) with respect to $f_i$ gives $a_{ij}=a_{ik}$.
Set $\alpha_i:=a_{ij}=a_{ik}.$
Then
\begin{gather*}
\cosh l_{ij}
=\sqrt{(1+\alpha_ie^{2f_i})(1+\alpha_je^{2f_j})}
+\eta_{ij}e^{f_i+f_j},\ \text{if}\ 1+\alpha_ie^{2f_i}>0,\ 1+\alpha_je^{2f_j}>0, \\
\cosh l_{ij}
=-\sqrt{(1+\alpha_ie^{2f_i})(1+\alpha_je^{2f_j})}
+\eta_{ij}e^{f_i+f_j},\ \text{if}\ 1+\alpha_ie^{2f_i}<0,\ 1+\alpha_je^{2f_j}<0,
\end{gather*}
with
\begin{equation*}
\tanh d_{ij}
=\frac{1}{\sinh l_{ij}}\frac{\partial}{\partial f_i}\cosh l_{ij}
=\frac{\alpha_ie^{2f_i}}{\sinh l_{ij}}
\sqrt{\frac{1+\alpha_je^{2f_j}}{1+\alpha_ie^{2f_i}}}
+\frac{\eta_{ij}e^{f_i+f_j}}{\sinh l_{ij}}.
\end{equation*}
This gives the discrete conformal structure $\mathbf{(b1)}$.
It is easy to check that the $\alpha_i$, $\eta_{ij}$ derived in any two adjacent triangles are equal.

\item[(II)] If $\frac{\partial H}{\partial f_i}\equiv 2$, then $\frac{\partial H}{\partial f_j}\equiv -2$ by (\ref{Eq: F8}).
This implies
\begin{equation}\label{Eq: F19}
H=2f_i-2f_j+c_{ij},
\end{equation}
where $c_{ij}$ is a constant. By the definition of $H$ in (\ref{defn of H}), we have $c_{ij}+c_{ji}=0$.
Then
\begin{equation*}
\frac{\cosh^2 d_{ij}}{\cosh^2 d_{ji}}
=e^H=c_2e^{2f_i-2f_j},
\end{equation*}
where $c_2=e^{c_{ij}}>0$.
Thus
$\frac{\cosh d_{ij}}{\cosh d_{ji}}=c_3e^{f_i-f_j}$ and
$\frac{\cosh d_{ji}}{\cosh d_{ij}}=c_4e^{f_j-f_i}$,
where $c_3=\sqrt{c_2}=e^{\frac{1}{2}c_{ij}}>0$ and $c_4=1/c_3=e^{-\frac{1}{2}c_{ij}}=e^{\frac{1}{2}c_{ji}}>0$.
By (\ref{Eq: F1}), we have
\begin{equation*}
\frac{\partial}{\partial f_i}\cosh l_{ij}
=\cosh l_{ij}-c_4e^{f_j-f_i},
\end{equation*}
which implies
$\cosh l_{ij}=\frac{1}{2}c_4 e^{f_j-f_i}+c_5(f_j)e^{f_i}$.
Similarly, by (\ref{Eq: F2}), we have
\begin{equation*}
\frac{\partial}{\partial f_j}\cosh l_{ij}
=\cosh l_{ij}-c_3e^{f_i-f_j},
\end{equation*}
which implies
$\cosh l_{ij}=\frac{1}{2}c_3 e^{f_i-f_j}+c_6(f_i)e^{f_j}$.
Hence,
\begin{equation*}
\frac{1}{2}c_4 e^{f_j-f_i}+c_5(f_j)e^{f_i}
=\frac{1}{2}c_3 e^{f_i-f_j}+c_6(f_i)e^{f_j},
\end{equation*}
which implies
$c_5(f_j)e^{-f_j}-\frac{1}{2}c_3 e^{-2f_j}
=c_6(f_i)e^{-f_i}-\frac{1}{2}c_4 e^{-2f_i}$.
Note that
$c_5(f_j)e^{-f_j}-\frac{1}{2}c_3 e^{-2f_j}$ depends only on $f_j$ and
$c_6(f_i)e^{-f_i}-\frac{1}{2}c_4 e^{-2f_i}$ depends only on $f_i$, there exists a constant
$\eta_{ij}$ such that $\eta_{ij}=c_5(f_j)e^{-f_j}-\frac{1}{2}c_3 e^{-2f_j}
=c_6(f_i)e^{-f_i}-\frac{1}{2}c_4 e^{-2f_i}$.
Then
$c_5(f_j)=\eta_{ij}e^{f_j}+\frac{1}{2}c_3 e^{-f_j}$,
$c_6(f_i)=\eta_{ij}e^{f_i}+\frac{1}{2}c_4 e^{-f_i}$.
Therefore,
\begin{equation*}
\cosh l_{ij}=\cosh (f_j-f_i-C_{ij})+\eta_{ij}e^{f_i+f_j},
\end{equation*}
where $C_{ij}=\log c_3=\frac{1}{2}c_{ij}$.
By (\ref{Eq: hyperbolic compatible condition}), we have
\begin{equation*}
\log\frac{\cosh^2 d_{ij}}{\cosh^2 d_{ji}}
+\log\frac{\cosh^2 d_{jk}}{\cosh^2 d_{kj}}
+\log\frac{\cosh^2 d_{ki}}{\cosh^2 d_{ik}}
=0.
\end{equation*}
Then $c_{ij}+c_{jk}+c_{ki}=0$ by (\ref{Eq: F19}).
Hence, $C_{ij}+C_{jk}+C_{ki}=0$ and $C_{rs}+C_{sr}=0$ for $\{r,s\}\subseteq\{i,j,k\}$.
This gives the discrete conformal structure $\mathbf{(b2)}$.
\end{description}
\qed

\section{Proof of spherical case in Theorem \ref{Thm: classification of DCS} $\mathbf{(C)}$}\label{section 3}
For the discrete conformal structure $\mathbf{(c1)}$, Glickenstein-Thomas \cite{GT} did not give a detailed proof, although the proof is similar.
For completeness, we provide the proof of the case $\mathbf{(C)}$ in details.
By the formula (\ref{Eq: spherical DCS}), we have
\begin{equation}\label{Eq: F10}
\frac{\partial}{\partial f_i}\cos l_{ij}
=-\sin l_{ij}\frac{\partial l_{ij}}{\partial f_i}
=-\sin l_{ij}\tan d_{ij}
=\cos l_{ij}-\frac{\cos d_{ji}}{\cos d_{ij}}.
\end{equation}
Similarly,
\begin{equation}\label{Eq: F11}
\frac{\partial}{\partial f_j}\cos l_{ij}
=\cos l_{ij}-\frac{\cos d_{ij}}{\cos d_{ji}}.
\end{equation}
Note that
\begin{equation*}
\frac{\partial}{\partial f_j}\tan d_{ij}
=\frac{\partial^2 l_{ij}}{\partial f_i\partial f_j}
=\frac{\partial}{\partial f_i}\tan d_{ji},
\end{equation*}
which is equivalent to
\begin{equation}\label{Eq: F12}
\frac{1}{\cos^2 d_{ij}}\frac{\partial d_{ij}}{\partial f_j}=\frac{1}{\cos^2 d_{ji}}\frac{\partial d_{ji}}{\partial f_i}.
\end{equation}
Set
\begin{equation}\label{defn of H spherical}
H=\log\frac{\cos^2 d_{ij}}{\cos^2 d_{ji}}.
\end{equation}
Then the formulas (\ref{Eq: spherical compatible condition}) and (\ref{Eq: variation}) give
\begin{equation}\label{Eq: F13}
\frac{\partial^2 H}{\partial f_i\partial f_j}=0.
\end{equation}
Direct calculations give
\begin{equation*}
\begin{aligned}
\frac{\partial H}{\partial f_i}
=&2\frac{\cos d_{ji}}{\cos d_{ij}}\frac{\partial}{\partial f_i}(\frac{\cos d_{ij}}{\cos d_{ji}})\\
=&2(-\tan d_{ij}\frac{\partial d_{ij}}{\partial f_i}+\tan d_{ji}\frac{\partial d_{ji}}{\partial f_i})\\
=&2\left[-\tan d_{ij}\frac{\partial d_{ij}}{\partial f_i}+\tan d_{ji}\bigg(\frac{\cos^2 d_{ji}}{\cos^2 d_{ij}}\bigg)\frac{\partial d_{ij}}{\partial f_j}\right],
\end{aligned}
\end{equation*}
where (\ref{Eq: F12}) is used in the last line.
Similarly,
\begin{equation*}
\begin{aligned}
\frac{\partial H}{\partial f_j}
=&2(-\tan d_{ij}\frac{\partial d_{ij}}{\partial f_j}+\tan d_{ji}\frac{\partial d_{ji}}{\partial f_j})\\
=&2\left[-\tan d_{ij}\bigg(\frac{\cos^2 d_{ij}}{\cos^2 d_{ji}}\bigg)\frac{\partial d_{ji}}{\partial f_i}+\tan d_{ji}\frac{\partial d_{ji}}{\partial f_j}\right].
\end{aligned}
\end{equation*}
Hence,
\begin{equation*}
\begin{aligned}
\bigg(\frac{\cos^2 d_{ij}}{\cos^2 d_{ji}}\frac{\partial}{\partial f_i}+\frac{\partial}{\partial f_j}\bigg)H
=&2\left[-\tan d_{ij}\bigg(\frac{\cos^2 d_{ij}}{\cos^2 d_{ji}}\bigg)\frac{\partial l_{ij}}{\partial f_i}+\tan d_{ji}\frac{\partial l_{ij}}{\partial f_j}\right]\\
=&2\left[-\tan^2 d_{ij}\bigg(\frac{\cos^2 d_{ij}}{\cos^2 d_{ji}}\bigg)+\tan^2 d_{ji}\right]\\
=&2\bigg(\frac{\cos^2 d_{ij}}{\cos^2 d_{ji}}-1\bigg),
\end{aligned}
\end{equation*}
which is equivalent to
\begin{equation}\label{Eq: F14}
\bigg(e^H\frac{\partial}{\partial f_i}+\frac{\partial}{\partial f_j}\bigg)H=2(e^H-1),
\end{equation}
and
\begin{equation}\label{Eq: F15}
\bigg(\frac{\partial}{\partial f_i}+e^{-H}\frac{\partial}{\partial f_j}\bigg)H=2(1-e^{-H}).
\end{equation}
Differentiating (\ref{Eq: F14}) with respect to $f_i$ and differentiating (\ref{Eq: F15}) with respect to $f_j$ give
\begin{equation}\label{Eq: partial H partial fi sphere}
\frac{\partial^2 H}{\partial f^2_i}+\bigg(\frac{\partial H}{\partial f_i}\bigg)^2=2\frac{\partial H}{\partial f_i},
\end{equation}
and
\begin{equation}\label{Eq: partial H partial fj sphere}
\frac{\partial^2 H}{\partial f^2_j}-\bigg(\frac{\partial H}{\partial f_j}\bigg)^2=2\frac{\partial H}{\partial f_j}.
\end{equation}
One can easily solve the ODE (\ref{Eq: partial H partial fi sphere}) to obtain that
\begin{equation*}
\frac{\partial H}{\partial f_i}\equiv2\ \quad \text{or}\ \quad
\frac{\partial H}{\partial f_i}=-2\frac{a_{ij}e^{2f_i}}{1-a_{ij}e^{2f_i}}
\end{equation*}
for some constants $a_{ij}$ by (\ref{Eq: F13}), and solve the ODE (\ref{Eq: partial H partial fj sphere}) to obtain that
\begin{equation*}
\frac{\partial H}{\partial f_j}\equiv -2\ \quad \text{or}\ \quad
\frac{\partial H}{\partial f_j}=2\frac{a_{ji}e^{2f_j}}{1-a_{ji}e^{2f_j}}
\end{equation*}
for some constants $a_{ji}$ by (\ref{Eq: F13}).
Note that the formula (\ref{Eq: F14}) can be rewritten as
\begin{equation}\label{Eq: F17}
(\frac{\partial H}{\partial f_i}-2)e^H
=-(\frac{\partial H}{\partial f_j}+2).
\end{equation}
\begin{description}
\item[(I)] If $\frac{\partial H}{\partial f_i}\not\equiv 2$, then $\frac{\partial H}{\partial f_j}\not\equiv 2$ by (\ref{Eq: F17}). Hence,
\begin{equation*}
\frac{\cos^2 d_{ij}}{\cos^2 d_{ji}}
=e^H=\frac{1}{\frac{\partial H}{\partial f_i}-2}(-\frac{\partial H}{\partial f_j}-2)=\frac{1-a_{ij}e^{2f_i}}{1-a_{ji}e^{2f_j}},
\end{equation*}
which implies
\begin{equation*}
\frac{\cos d_{ij}}{\cos d_{ji}}
=\pm\sqrt{\frac{1-a_{ij}e^{2f_i}}{1-a_{ji}e^{2f_j}}}.
\end{equation*}
\begin{description}
\item[(i)]
If $\frac{\cos d_{ij}}{\cos d_{ji}}
=\sqrt{\frac{1-a_{ij}e^{2f_i}}{1-a_{ji}e^{2f_j}}}$,
then by (\ref{Eq: F10}) and (\ref{Eq: F11}), one can obtain the following system
\begin{equation*}
\begin{aligned}
&\frac{\partial}{\partial f_i}\cos l_{ij}
=\cos l_{ij}
-\sqrt{\frac{1-a_{ji}e^{2f_j}}{1-a_{ij}e^{2f_i}}},\\
&\frac{\partial}{\partial f_j}\cos l_{ij}
=\cos l_{ij}
-\sqrt{\frac{1-a_{ij}e^{2f_i}}{1-a_{ji}e^{2f_j}}}.
\end{aligned}
\end{equation*}
This implies
\begin{equation*}
\cos l_{ij}
=\sqrt{(1-a_{ij}e^{2f_i})(1-a_{ji}e^{2f_j})}
-\eta_{ij}e^{f_i+f_j},\ \text{if}\ 1-a_{ij}e^{2f_i}>0,\ 1-a_{ji}e^{2f_j}>0,
\end{equation*}
\begin{equation*}
\cos l_{ij}
=-\sqrt{(1-a_{ij}e^{2f_i})(1-a_{ji}e^{2f_j})}
-\eta_{ij}e^{f_i+f_j},\ \text{if}\ 1-a_{ij}e^{2f_i}<0,\ 1-a_{ji}e^{2f_j}<0,
\end{equation*}
where $\eta_{ij}$ is constant.
By (\ref{Eq: spherical compatible condition}) and (\ref{Eq: F13}), we have
\begin{equation}\label{Eq: F16}
\log\frac{\cos d_{ij}}{\cos d_{ji}}
+\log\frac{\cos d_{ki}}{\cos d_{ik}}
=\log\frac{\cos d_{kj}}{\cos d_{jk}}.
\end{equation}
Note that $\log\frac{\cos d_{kj}}{\cos d_{jk}}$ is independent of $f_i$ by (\ref{Eq: variation}), so differentiating (\ref{Eq: F16}) with respect to $f_i$ gives $a_{ij}=a_{ik}$.
Therefore, one can set $\alpha_i:=a_{ij}=a_{ik}$,
and then
\begin{equation*}
\cos l_{ij}
=\sqrt{(1-\alpha_ie^{2f_i})(1-\alpha_je^{2f_j})}
-\eta_{ij}e^{f_i+f_j},\ \text{if}\ 1-\alpha_ie^{2f_i}>0,\ 1-\alpha_je^{2f_j}>0,
\end{equation*}
\begin{equation*}
\cos l_{ij}
=-\sqrt{(1-\alpha_ie^{2f_i})(1-\alpha_je^{2f_j})}
-\eta_{ij}e^{f_i+f_j},\ \text{if}\ 1-\alpha_ie^{2f_i}<0,\ 1-\alpha_je^{2f_j}<0,
\end{equation*}
with
\begin{equation*}
\tan d_{ij}
=-\frac{1}{\sin l_{ij}}\frac{\partial}{\partial f_i}\cos l_{ij}
=\frac{\alpha_ie^{2f_i}}{\sin l_{ij}}\sqrt{\frac{1-\alpha_je^{2f_j}}{1-\alpha_ie^{2f_i}}}
+\frac{\eta_{ij}e^{f_i+f_j}}{\sin l_{ij}}.
\end{equation*}
It is easy to check that the $\alpha_i$ and $\eta_{ij}$ derived in any two adjacent spherical triangles must be equal.
This gives the discrete conformal structure $\mathbf{(c1)}$.

\item[(ii)]
If $\frac{\cos d_{ij}}{\cos d_{ji}}
=-\sqrt{\frac{1-a_{ij}e^{2f_i}}{1-a_{ji}e^{2f_j}}}$,
then one can obtain the following system
\begin{equation*}
\begin{aligned}
&\frac{\partial}{\partial f_i}\cos l_{ij}
=\cos l_{ij}
+\sqrt{\frac{1-a_{ji}e^{2f_j}}{1-a_{ij}e^{2f_i}}},\\
&\frac{\partial}{\partial f_j}\cos l_{ij}
=\cos l_{ij}
+\sqrt{\frac{1-a_{ij}e^{2f_i}}{1-a_{ji}e^{2f_j}}}.
\end{aligned}
\end{equation*}
This implies
\begin{equation*}
\cos l_{ij}
=-\sqrt{(1-a_{ij}e^{2f_i})(1-a_{ji}e^{2f_j})}
+\eta_{ij}e^{f_i+f_j},\ \text{if}\ 1-a_{ij}e^{2f_i}>0,\ 1-a_{ji}e^{2f_j}>0,
\end{equation*}
\begin{equation*}
\cos l_{ij}
=\sqrt{(1-a_{ij}e^{2f_i})(1-a_{ji}e^{2f_j})}
+\eta_{ij}e^{f_i+f_j},\ \text{if}\ 1-a_{ij}e^{2f_i}<0,\ 1-a_{ji}e^{2f_j}<0,
\end{equation*}
where $\eta_{ij}$ is a constant.
Similarly, one can also set $\alpha_i:=a_{ij}=a_{ik}$,
and then
\begin{equation*}
\cos l_{ij}
=-\sqrt{(1-\alpha_ie^{2f_i})(1-\alpha_je^{2f_j})}
+\eta_{ij}e^{f_i+f_j},\ \text{if}\ 1-\alpha_ie^{2f_i}>0,\ 1-\alpha_je^{2f_j}>0,
\end{equation*}
\begin{equation*}
\cos l_{ij}
=\sqrt{(1-\alpha_ie^{2f_i})(1-\alpha_je^{2f_j})}
+\eta_{ij}e^{f_i+f_j},\ \text{if}\ 1-\alpha_ie^{2f_i}<0,\ 1-\alpha_je^{2f_j}<0,
\end{equation*}
with
\begin{equation*}
\tan d_{ij}
=-\frac{\alpha_ie^{2f_i}}{\sinh l_{ij}}\sqrt{\frac{1-\alpha_je^{2f_j}}{1-\alpha_ie^{2f_i}}}
-\frac{\eta_{ij}e^{f_i+f_j}}{\sin l_{ij}}.
\end{equation*}
It is easy to check that the $\alpha_i$ and $\eta_{ij}$ derived in any two adjacent spherical triangles must be equal.
This gives the discrete conformal structure $\mathbf{(c3)}$.
\end{description}

\item[(II)]
If $\frac{\partial H}{\partial f_i}\equiv 2$, then $\frac{\partial H}{\partial f_j}\equiv -2$ by (\ref{Eq: F17}).
This implies
\begin{equation}\label{Eq: F20}
H=2f_i-2f_j+c_{ij},
\end{equation}
where $c_{ij}$ is a constant. By the definition of $H$ in (\ref{defn of H spherical}), we have $c_{ij}+c_{ji}=0$.
Then
\begin{equation*}
\frac{\cos^2 d_{ij}}{\cos^2 d_{ji}}
=e^H=c_2e^{2f_i-2f_j},
\end{equation*}
where $c_2=e^{c_{ij}}>0$.
\begin{description}
\item[(i)] If
$\frac{\cos d_{ij}}{\cos d_{ji}}=c_3e^{f_i-f_j}$ and
$\frac{\cos d_{ji}}{\cos d_{ij}}=c_4e^{f_j-f_i}$,
where $c_3=\sqrt{c_2}=e^{\frac{1}{2}c_{ij}}>0$ and $c_4=1/c_3=e^{\frac{1}{2}c_{ji}}>0$.
Then by (\ref{Eq: F10}), we have
\begin{equation*}
\frac{\partial}{\partial f_i}\cos l_{ij}
=\cos l_{ij}-c_4e^{f_j-f_i},
\end{equation*}
which implies
$\cos l_{ij}=\frac{1}{2}c_4 e^{f_j-f_i}+c_5(f_j)e^{f_i}$.
Similarly, by (\ref{Eq: F11}), we have
\begin{equation*}
\frac{\partial}{\partial f_j}\cos l_{ij}
=\cos l_{ij}-c_3e^{f_i-f_j},
\end{equation*}
which implies
$\cos l_{ij}=\frac{1}{2}c_3 e^{f_i-f_j}+c_6(f_i)e^{f_j}$.
Hence,
\begin{equation*}
\frac{1}{2}c_4 e^{f_j-f_i}+c_5(f_j)e^{f_i}
=\frac{1}{2}c_3 e^{f_i-f_j}+c_6(f_i)e^{f_j},
\end{equation*}
which implies
$c_5(f_j)e^{-f_j}-\frac{1}{2}c_3 e^{-2f_j}
=c_6(f_i)e^{-f_i}-\frac{1}{2}c_4 e^{-2f_i}:=-\eta_{ij}$.
Then
$c_5(f_j)=-\eta_{ij}e^{f_j}+\frac{1}{2}c_3 e^{-f_j}$,
$c_6(f_i)=-\eta_{ij}e^{f_i}+\frac{1}{2}c_4 e^{-f_i}$.
Therefore,
\begin{equation*}
\cos l_{ij}=\cosh (f_j-f_i-C_{ij})-\eta_{ij}e^{f_i+f_j},
\end{equation*}
where $C_{ij}=c_3=\frac{1}{2}c_{ij}$.
By (\ref{Eq: spherical compatible condition}), we have
\begin{equation*}
\log\frac{\cos^2 d_{ij}}{\cos^2 d_{ji}}
+\log\frac{\cos^2 d_{jk}}{\cos^2 d_{kj}}
+\log\frac{\cos^2 d_{ki}}{\cos^2 d_{ik}}
=0.
\end{equation*}
Then $c_{ij}+c_{jk}+c_{ki}=0$ by (\ref{Eq: F20}).
Hence, $C_{ij}+C_{jk}+C_{ki}=0$ and $C_{rs}+C_{sr}=0$ for $\{r,s\}\subseteq\{i,j,k\}$.
This gives the discrete conformal structure $\mathbf{(c2)}$.
\item[(ii)]
If $\frac{\cos d_{ij}}{\cos d_{ji}}=-c_3e^{f_i-f_j}$ and $\frac{\cos d_{ji}}{\cos d_{ij}}=-c_4e^{f_j-f_i}$.
Then by (\ref{Eq: F10}), we have
\begin{equation*}
\frac{\partial}{\partial f_i}\cos l_{ij}
=\cos l_{ij}+c_4e^{f_j-f_i},
\end{equation*}
which implies
$\cos l_{ij}=-\frac{1}{2}c_4 e^{f_j-f_i}+c_5(f_j)e^{f_i}$.
Similarly, by (\ref{Eq: F11}), we have
\begin{equation*}
\frac{\partial}{\partial f_j}\cos l_{ij}
=\cos l_{ij}+c_3e^{f_i-f_j},
\end{equation*}
which implies
$\cos l_{ij}=-\frac{1}{2}c_3 e^{f_i-f_j}+c_{6}(f_i)e^{f_j}$.
Hence,
\begin{equation*}
-\frac{1}{2}c_4 e^{f_j-f_i}+c_5(f_j)e^{f_i}
=-\frac{1}{2}c_3 e^{f_i-f_j}+c_6(f_i)e^{f_j},
\end{equation*}
which implies
$c_5(f_j)e^{-f_j}+\frac{1}{2}c_3 e^{-2f_j}
=c_6(f_i)e^{-f_i}+\frac{1}{2}c_4 e^{-2f_i}:=\eta_{ij}$.
Then
$c_5(f_j)=\eta_{ij}e^{f_j}-\frac{1}{2}c_3 e^{-f_j}$,
$c_6(f_i)=\eta_{ij}e^{f_i}-\frac{1}{2}c_4 e^{-f_i}$.
Therefore
\begin{equation*}
\cos l_{ij}=-\cosh (f_j-f_i-C_{ij})+\eta_{ij}e^{f_i+f_j},
\end{equation*}
where $C_{ij}=\log c_3=\frac{1}{2}c_{ij}$.
Similarly, we have
$C_{ij}+C_{jk}+C_{ki}=0$ and $C_{rs}+C_{sr}=0$ for $\{r,s\}\subseteq\{i,j,k\}$.
This gives the discrete conformal structure $\mathbf{(c4)}$.
\end{description}

\end{description}
\qed

\section{Relationships with 3-dimensional hyperbolic geometry}\label{section 4}
In this section, we study the relationships between the discrete conformal structures on closed surfaces in Theorem \ref{Thm: classification of DCS} and 3-dimensional hyperbolic geometry.

The relationships between the discrete conformal structures on closed triangulated surfaces and 3-dimensional hyperbolic geometry were first observed by Bobenko-Pinkall-Springborn \cite{BPS} in the case of Luo's vertex scaling of piecewise linear metrics.
Suppose $Ov_iv_jv_k$ is an ideal hyperbolic tetrahedron in $\mathbb{H}^3$ with each ideal vertex attached with a horosphere.
Bobenko-Pinkall-Springborn \cite{BPS} found that Luo's constuction of Euclidean triangle via vertex scaling corresponds exactly to the Euclidean triangle given by the intersection of $Ov_iv_jv_k$ and the horosphere at the ideal vertex $O$, if the generalized edge lengths of the decorated ideal hyperbolic tetrahedron $Ov_iv_jv_k$ are properly assigned.
Based on this observation, Bobenko-Pinkall-Springborn \cite{BPS} further introduced the vertex scaling for piecewise hyperbolic (spherical respectively) metrics by perturbing the ideal vertex $O$ of the ideal hyperbolic tetrahedron
$Ov_iv_jv_k$ to hyper-ideal (hyperbolic respectively) while keeping the other vertices ideal.
Motivated by Bobenko-Pinkall-Springborn's observations \cite{BPS}, Zhang-Guo-Zeng-Luo-Yau-Gu \cite{ZGZLYG} further constructed all 18 types of discrete conformal structures on closed surfaces in $\mathbb{G}$ background geometries by perturbing the ideal vertices of the ideal hyperbolic tetrahedron $Ov_iv_jv_k$ to be hyperbolic, ideal or hyper-ideal.
Furthermore, Zhang-Guo-Zeng-Luo-Yau-Gu \cite{ZGZLYG} introduced the following definition of discrete conformal structure on surfaces.
\begin{definition}[\cite{ZGZLYG}, Definition 3.9]\label{Def: ZGZLYG}
Suppose $(M,\mathcal{T})$ is a triangulated surface with $\mathbb{G}$ background geometry.
Let $\varepsilon: V\rightarrow \{-1,0,1\}$ and $\zeta\in E^*$ be two weights defined on the vertices and edges respectively with $\zeta_{ij}=\zeta_{ji}$.
A discrete conformal structure on $(M,\mathcal{T})$ is a map $u\in V^*$ such that the edge length $l_{ij}$ for the edge $\{ij\}\in E$ is given by
\begin{equation}\label{Eq: Euclidean edge}
l_{ij}^2=2\zeta_{ij}e^{u_i+u_j}+\varepsilon_ie^{2u_i}+\varepsilon_je^{2u_j}
\quad \text{if $\mathbb{G}=\mathbb{E}^2$},
\end{equation}
\begin{equation}\label{Eq: hyperbolic edge}
\cosh l_{ij}=\frac{4\zeta_{ij}e^{u_i+u_j}
+(1+\varepsilon_ie^{2u_i})(1+\varepsilon_je^{2u_j})}
{(1-\varepsilon_ie^{2u_i})(1-\varepsilon_je^{2u_j})}
\quad \text{if $\mathbb{G}=\mathbb{H}^2$},
\end{equation}
\begin{equation}\label{Eq: spherical edge}
\cos l_{ij}=\frac{-4\zeta_{ij}e^{u_i+u_j}
+(1-\varepsilon_ie^{2u_i})(1-\varepsilon_je^{2u_j})}
{(1+\varepsilon_ie^{2u_i})(1+\varepsilon_je^{2u_j})}
\quad \text{if $\mathbb{G}=\mathbb{S}^2$}.
\end{equation}
\end{definition}

By a change of variables, Glickenstein-Thomas' classification of discrete conformal structures in \cite{GT} is equivalent to that in Definition \ref{Def: ZGZLYG}.
By Remark \ref{remark: 2}, we set $\alpha_i\in \{-1,0,1\}$.
Hence, it is obvious that the discrete conformal structure $\mathbf{(A)}$ in Theorem \ref{Thm: classification of DCS} is equivalent to (\ref{Eq: Euclidean edge}) in Definition \ref{Def: ZGZLYG}.
Moreover, on a triangulated surface $(M,\mathcal{T})$ with $\mathbb{H}^2$ background geometry, if we set
\begin{eqnarray*}
e^{f_i}=
\begin{cases}
r_i  &{\text{if $\alpha_i=0$}},\\
\sinh r_i &{\text{if $\alpha_i=1$}},\\
\tanh r_i &{\text{if $\alpha_i=-1$}}
\end{cases}
\end{eqnarray*}
for the discrete conformal structure (\ref{Eq: F22}) in Theorem \ref{Thm: classification of DCS} and set
\begin{eqnarray*}
e^{u_i}=
\begin{cases}
r_i  &{\text{if $\varepsilon_i=0$}},\\
\tanh\frac{r_i}{2} &{\text{if $\varepsilon_i=\pm1$}}
\end{cases}
\end{eqnarray*}
for the discrete conformal structure (\ref{Eq: hyperbolic edge}) in Definition \ref{Def: ZGZLYG},
then it is easy to check that the discrete conformal structure (\ref{Eq: F22}) in Theorem \ref{Thm: classification of DCS} is equivalent to (\ref{Eq: hyperbolic edge}) in Definition \ref{Def: ZGZLYG}.
Similarly, on a triangulated surface $(M,\mathcal{T})$ with $\mathbb{S}^2$ background geometry, if we set
\begin{eqnarray*}
e^{f_i}=
\begin{cases}
r_i  &{\text{if $\alpha_i=0$}},\\
\sin r_i &{\text{if $\alpha_i=1$}},\\
\tan r_i &{\text{if $\alpha_i=-1$}}
\end{cases}
\end{eqnarray*}
for the discrete conformal structure (\ref{Eq: F32}) in Theorem \ref{Thm: classification of DCS} and set
\begin{eqnarray*}
e^{u_i}=
\begin{cases}
r_i  &{\text{if $\varepsilon_i=0$}},\\
\tan\frac{r_i}{2} &{\text{if $\varepsilon_i=\pm1$}}
\end{cases}
\end{eqnarray*}
for the discrete conformal structure (\ref{Eq: spherical edge}) in Definition \ref{Def: ZGZLYG},
then it is easy to check that the discrete conformal structure (\ref{Eq: F32}) in Theorem \ref{Thm: classification of DCS} is equivalent to (\ref{Eq: spherical edge}) in Definition \ref{Def: ZGZLYG}.
For example, if we set $\alpha\equiv 1$ and $e^f=\sin r$, then the discrete conformal structure (\ref{Eq: F32}) can be reduced to
\begin{equation*}
\cos l_{ij}
=\cos r_i\cos r_j-\eta_{ij}\sin r_i\sin r_j,
\end{equation*}
which is equivalent to (\ref{Eq: spherical edge}) for $\varepsilon\equiv 1$ and $e^u=\tan \frac{r}{2}$.
Note that the discrete conformal structure (\ref{Eq: spherical edge}) for $\varepsilon\equiv 1$ corresponds to  the cosine law for a generalized hyperbolic triangle with two hyper-ideal vertices and one hyperbolic vertex and
$u$ represents the edge length of this generalized hyperbolic triangle.
Therefore, the discrete conformal structure (\ref{Eq: F32}) in Theorem \ref{Thm: classification of DCS} takes the edge lengths as its variables.

While the discrete conformal structure (\ref{Eq: F34}) in Theorem \ref{Thm: classification of DCS} takes the generalized angles as its variables, which can be taken as the dual form of (\ref{Eq: F32}).
Set
\begin{eqnarray*}
e^{f_i}=
\begin{cases}
r_i  &{\text{if $\alpha_i=0$}},\\
\sin r_i &{\text{if $\alpha_i=1$}},\\
\sinh r_i &{\text{if $\alpha_i=-1$}}.
\end{cases}
\end{eqnarray*}

In the case that $\alpha\equiv0$, the discrete conformal structure (\ref{Eq: F34}) can be reduced to
\begin{equation*}
\cos l_{ij}=-1+\eta_{ij}r_ir_j,
\end{equation*}
which corresponds to the cosine law for a generalized hyperbolic triangle with two ideal vertices and one hyperbolic vertex.

In the case that $\alpha\equiv1$, the discrete conformal structure (\ref{Eq: F34}) can be reduced to
\begin{equation*}
\cos l_{ij}=-\cos r_i\cos r_j+\eta_{ij}\sin r_i\sin r_j,
\end{equation*}
which corresponds to the cosine law for a generalized hyperbolic triangle with three hyperbolic vertices if $\eta_{ij}>1$.

In the case that $\alpha\equiv-1$, the discrete conformal structure (\ref{Eq: F34}) can be reduced to
\begin{equation*}
\cos l_{ij}=-\cosh r_i\cosh r_j+\eta_{ij}\sinh r_i\sinh r_j,
\end{equation*}
which corresponds to the cosine law for a generalized hyperbolic triangle with two hyper-ideal vertices and one hyperbolic vertex if $\eta_{ij}>1$.

In the case that $\alpha_i=0$ and $\alpha_j=1$, the discrete conformal structure (\ref{Eq: F34}) can be reduced to
\begin{equation*}
\cos l_{ij}=-\cos r_j+\eta_{ij} r_i\sin r_j,
\end{equation*}
which corresponds to the cosine law for a generalized hyperbolic triangle with two hyperbolic vertices and one ideal vertex if $\eta_{ij}>0$.

In the case that $\alpha_i=0$ and $\alpha_j=-1$, the discrete conformal structure (\ref{Eq: F34}) can be reduced to
\begin{equation*}
\cos l_{ij}=-\cosh r_j+\eta_{ij} r_i\sinh r_j,
\end{equation*}
which corresponds to the cosine law for a generalized hyperbolic triangle with a hyper-ideal vertex, a ideal vertex and a hyperbolic vertex if $\eta_{ij}>0$.

In the case that $\alpha_i=1$ and $\alpha_j=-1$, the discrete conformal structure (\ref{Eq: F34}) can be reduced to
\begin{equation*}
\cos l_{ij}=-\cos r_i\cosh r_j+\eta_{ij}\sin r_i\sinh r_j,
\end{equation*}
which corresponds to the cosine law for a generalized hyperbolic triangle with two hyperbolic vertices and one hyper-ideal vertex if $\eta_{ij}>0$.

In the six cases above, $r$ represents the generalized angle of a generalized hyperbolic triangle.

For the discrete conformal structure (\ref{Eq: F23}), we can set $\alpha\equiv 1$ and
$e^{f_i}=\cosh r_i$, then 
\begin{equation*}
\cosh l_{ij}=-\sinh r_i\sinh r_j+\eta_{ij}\cosh r_i\cosh r_j,
\end{equation*}
which corresponds to the cosine law for a generalized hyperbolic triangle with two hyperbolic vertices and one hyper-ideal vertex with 
$r$ representing the length.
For the discrete conformal structure $\mathbf{(b2)}$ with $C\equiv 0$ in Theorem \ref{Thm: classification of DCS},
it is easy to check that it corresponds to the cosine law for a generalized hyperbolic triangle with two hyperbolic vertices and one ideal vertex.
However, the relationships between the hyperbolic geometry and the discrete conformal structures (\ref{Eq: F33}), $\mathbf{(c2)}$, (\ref{Eq: F35})and $\mathbf{(c4)}$ are not clear.

By the arguments above, the discrete conformal structures in Theorem \ref{Thm: classification of DCS} contain some of the generalized circle packings introduced by Guo-Luo \cite{GL2} as special cases.
Please refer to Guo-Luo \cite{GL2} for more details.

\end{document}